\documentclass[12pt]{article}
\usepackage{amsmath,amssymb,amsthm}
\usepackage{vmargin}
\usepackage{url}

\setpapersize{USletter}
\setmargnohfrb{1.0in}{1.0in}{1.0in}{1.0in}

\usepackage{pxfonts}

\newtheorem{Th}{Theorem}

\newcommand{\lcm}{\mathop{\mathrm{lcm}}}
\newcommand{\beq}{\begin{equation}}
\newcommand{\eeq}{\end{equation}}
\newcommand{\bea}{\begin{array}}
\newcommand{\eea}{\end{array}}
\newcommand{\eqdef}{\mathbin{\mathop=\limits^{\rm def}}}

\title{On the intersections of Fibonacci, Pell, and Lucas numbers}
\author{Max A. Alekseyev\thanks{Department of Computer Science and Engineering,
University of South Carolina.\newline{Email: maxal@cse.sc.edu}} }

\begin{document}
\maketitle\thispagestyle{empty}

\begin{abstract}
We describe how to compute the intersection of two Lucas sequences of the forms 
\(\{ U_n(P,\pm 1) \}_{n=0}^{\infty}\) or $\{ V_n(P,\pm 1) \}_{n=0}^{\infty}$ with 
$P\in\mathbb{Z}$ that includes sequences of Fibonacci, Pell, Lucas, and Lucas-Pell numbers.
We prove that such an intersection is finite except for the case $U_n(1,-1)$ and $U_n(3,1)$ and the case 
of two $V$-sequences when the product of their discriminants is a perfect square. Moreover, 
the intersection in these cases also forms a Lucas sequence.
Our approach relies on solving homogeneous quadratic Diophantine equations and Thue equations.
In particular, we prove that $0, 1, 2,$ and $5$ are the only numbers that are both Fibonacci and Pell, 
and list similar results for many other pairs of Lucas sequences.
We further extend our results to Lucas sequences with arbitrary initial terms.
\end{abstract}

In 1962 Stein~\cite{Stein1962} using elementary arguments proved that if two Fibonacci sequences (with different initial terms) 
share three terms then they are identical from some terms on. This result was later generalized by Revuz~\cite{Revuz1974} 
to arbitrary sequences satisfying the same second order linear recurrence. A similar problem is to determine the intersection 
of two sequences satisfying distinct linear recurrences. 
Mignotte, 1978~\cite{Mignotte1978} proved that the intersection of two sequences is finite 
unless the roots of their characteristic polynomials are multiplicatively dependent. 
M\'aty\'as, 1981~\cite{Matyas1981} gave a criterion for determining whether two second order linear sequences have nonempty intersection.
Posing the problem as a system of Pellian equations, 
Pinch, 1988~\cite{Pinch1988} and Tzanakis, 2002~\cite{Tzanakis02} proposed computational methods for solving such systems.

In the current paper we focus on the computational matters of finding the intersection and provide some explicit results for a class of Lucas sequences.
While the general structure of such intersections is well known~\cite{Mignotte1978,Matyas1981,Pinch1988,Tzanakis02}, 
our approach (outlined in the proof of Theorem~\ref{Txyz}) is different from previous studies (see \cite{Tzanakis02} for a review of various computational methods). 
In short, we reduce the problem to a finite number of Thue equations, each of which has a finite number of solutions. 
In practice the intersection of particular Lucas sequences can be computed with PARI/GP computer algebra system~\cite{PARI}
providing a functionality for solving Thue equations, based on Bilu and Hanrot's improvement~\cite{Bilu1996}
of Tzanakis and de Weger's method~\cite{Tzanakis1989}.
We also characterize the cases when the intersection is infinite and show that it also forms a Lucas sequence.

While the current paper was under review, our attention was drawn to 
the paper of Szalay, 2007~\cite{Szalay2007}, which employs a
similar approach for solving systems of Pellian equations. 
We therefore feel obliged to underline the differences between 
\cite{Szalay2007} and our paper. 
In Theorem~\ref{Txyz} we consider a general system of two homogeneous 
quadratic equations in two indeterminants and 
formulate conditions under which such a system has a finite number of integer solutions.
The systems of Pellian equations considered in \cite{Szalay2007} 
represent a particular case of the system from our Theorem~\ref{Txyz},
while the algorithm and the corresponding proof of solutions' finiteness in \cite{Szalay2007} are flawed.\footnote{In \cite{Szalay2007} 
it is incorrectly claimed that the constructed quartic equation $(2.5)$ represents a Thue equation.
We found a counterexample with $(a_1,b_1,c_1)=(5,-1,-1)$ and $(a_2,b_2,c_2)=(20,-1,1)$ in $(2.1)-(2.2)$
which yields $(a,b,c)=(25,-1,-1)$ in $(2.4)$. From 
a basic solution $(X_0,Y_0,Z_0)=(1,0,5)$ it further leads to 
the reducible quartic polynomial 
$T(s,r) = 3125 s^4 - 2250 s^2 r^2 + 5 r^4 = 5 (25 s^2 - 20 sr - r^2) (25 s^2 + 20 sr - r^2)$ in the l.h.s. of $(2.5)$.}
In \cite{Szalay2007} it is used the result from \cite{Luca2006} 
that the denominator 
in the general solution to a Diophantine equation $Ax^2+By^2+Cz^2=0$, 
derived from a particular solution
$(x_0,y_0,z_0)$ with $z_0\ne 0$, divides $2A^2BCz_0^3$.
In Theorem~\ref{ThABC} we prove a stronger result
that this denominator actually divides $2\lcm(A,B) C z_0^2$,
leading to a smaller number of Thue equations to solve.
This result also corrects an error in Corollary 6.3.8 of \cite{Cohen07}.
We further focus on multiple issues specific to 
the intersections of Lucas sequences which are out of scope of \cite{Szalay2007}.
Said that, we refer to \cite{Szalay2007} for 
many examples of systems of Pellian equations and 
a review of their applications and related literature.

The paper is organized as follows. In Section~\ref{SecBasic} we give definition and basic properties of Lucas sequences.
We develop an algorithmic approach to solving systems of quadratic Diophantine equations in Section~\ref{SecQuad} 
and apply it to the problem of finding intersections of Lucas sequences in Section~\ref{SecInter}. 
Some numerical results are given and a number of open questions are posed in Section~\ref{SecComp}.
Finally, we discuss generalizations of our approach to Lucas sequences with arbitrary initial terms in Section~\ref{SecArbit}.

\section{Basic properties of Lucas sequences}\label{SecBasic}

The pair of Lucas sequences $U_n(P,Q)$ and $V_n(P,Q)$ are defined by the same linear recurrent relation with the coefficient $P,Q\in\mathbb{Z}$ but
different initial terms:
$$
\begin{array}{lll} 
U_0(P,Q)=0, & U_1(P,Q)=1, & U_{n+1}(P,Q) = P\cdot U_n(P,Q) - Q\cdot U_{n-1}(P,Q),\, n\geq 1;\\
V_0(P,Q)=2, & V_1(P,Q)=P, & V_{n+1}(P,Q) = P\cdot V_n(P,Q) - Q\cdot V_{n-1}(P,Q),\, n\geq 1.
\end{array}
$$

Some Lucas sequences have their own names:
\begin{center} \begin{tabular}{|l|l|l|}
\hline
Sequence & Name & Initial terms\\
\hline\hline
$U_n(1,-1)$ & Fibonacci numbers & $0, 1, 1, 2, 3, 5, 8, 13, 21, \dots$\\
\hline
$V_n(1,-1)$ & Lucas numbers & $2, 1, 3, 4, 7, 11, 18, 29, \dots$\\
\hline
$U_n(2,-1)$ & Pell numbers &  $0, 1, 2, 5, 12, 29, 70, 169, \dots$\\
\hline
$V_n(2,-1)$ & Pell-Lucas numbers &  $2, 2, 6, 14, 34, 82, 198, \dots$\\
\hline
\end{tabular}
\end{center}
Other examples include Jacobsthal numbers $U_n(1,-2)$, Mersenne numbers $U_n(3,2)$ etc.

In the current paper we focus on the case of $Q=1$ or $Q=-1$. 
We also notice that
$$U_n(-P,Q) = (-1)^{n+1} U_n(P,Q)\quad\text{and}\quad V_n(-P,Q) = (-1)^n V_n(P,Q)$$
and restrict our attention to Lucas sequences with $P\geq 0$.
Similarly, we exclude from consideration 
terms with negative indices, noticing that they may deviate only in signs:
\begin{center}
\begin{tabular}{ll}
$U_{-n}(P,1) = -U_n(P,1),$ & $V_{-n}(P,1)=V_n(P,1),$\\
$U_{-n}(P,-1) = (-1)^{n+1} U_n(P,-1),$ & $V_{-n}(P,-1)=(-1)^{n+1} V_n(P,1).$
\end{tabular}
\end{center}
From further consideration we also exclude the following degenerate cases:
\begin{center} \begin{tabular}{|l|l|l|}
\hline
$(P,Q)$ & $U_n(P,Q)$ & $V_n(P,Q)$ \\
\hline\hline
$(0,1)$ & $U_{2m}=0$ & $V_{2m}=2\cdot (-1)^m$ \\
        & $U_{2m+1}=(-1)^m$ & $V_{2m+1}=0$ \\
\hline
$(0,-1)$ & $U_{2m}=0$ & $V_{2m}=2$ \\
         & $U_{2m+1}=(-1)^m$ & $V_{2m+1}=0$ \\
\hline
$(1,1)$ & $U_{3m}=0$ & $V_{3m}=2$ \\
        & $U_{3m+1}=(-1)^m$ & $V_{3m+1}=(-1)^m$ \\
        & $U_{3m+2}=(-1)^m$ & $V_{3m+2}=(-1)^{m+1}$ \\
\hline
$(2,1)$ & $U_m = m$ & $V_m = 2$ \\
\hline
\end{tabular}
\end{center}
It is easy to see that under the described restrictions both $U_n(P,Q)$ and $V_n(P,Q)$ are positive for $n>0$.

The characteristic polynomial of Lucas sequences $\{ U_n(P,Q) \}$ and $\{ V_n(P,Q) \}$ is $\lambda^2 - P\lambda + Q$ with the discriminant $D=P^2 - 4Q$.
For non-degenerate sequences, the discriminant $D$ is a positive non-square integer.

Let $\alpha=\frac{P+\sqrt{D}}{2}$ and $\beta=\frac{P-\sqrt{D}}{2}$ be the roots of the characteristic polynomial, then the following explicit (Binet-type) formulas
take place
$$
\begin{array}{l}
Q = \alpha\beta \\
D = (\alpha - \beta)^2\\
U_n(P,Q) = \frac{\alpha^n - \beta^n}{\alpha - \beta} \\
V_n(P,Q) = \alpha^n + \beta^n
\end{array}
$$
In particular, these formulas imply that
\begin{equation}\label{MainEq}
V_n(P,Q)^2 - D\cdot U_n(P,Q)^2 = 4 Q^n.
\end{equation}
For $|Q|=1$, it means that the pairs $(V_n(P,Q),U_n(P,Q))$ form 
solutions to the equation:\footnote{Here and everywhere below $\pm$ in the r.h.s. of an equation means that we accept both signs as solutions.} 
\begin{equation}\label{EqXDY}
x^2 - D y^2 = \pm 4.
\end{equation}
The converse statement can be used to prove that given positive integers 
belong to $\{ V_n(P,Q) \}$ or $\{ U_n(P,Q) \}$ respectively 
(stated without a proof for Fibonacci numbers in \cite{Posamentier07}):

\begin{Th}\label{Trecog} 
Let $P$, $Q$ be integers such that $P>0$, $|Q|=1$, $(P,Q)\ne(3,1)$, 
and $D=P^2 - 4Q > 0$.
If positive integers $u$ and $v$ are such that 
$$v^2 - D u^2 = \pm 4,$$
then
$$u\in \{ U_n(P,Q) \}\qquad\text{and}\qquad v\in \{ V_n(P,Q) \}.$$
\end{Th}

\begin{proof} 
We notice that $D$ is congruent to $0$ or $1$ modulo $4$, implying 
that the general solution to equation \eqref{EqXDY} has the form $(x,y)=(x_k,y_k)$, where
$$\frac{x_k + y_k\sqrt{D}}2 = \pm \left(\frac{x_1 + y_1\sqrt{D}}2\right)^k,\qquad k\in\mathbb{Z},$$
and $(x_1,y_1)$ is a solution with the smallest positive $y_1$ 
(see Proposition 6.3.16 in \cite{Cohen07}). Since there is a solution $(x,y) = (P,1)$ with the smallest possible positive 
$y = 1$, we can take $(x_1,y_1) = (P,1)$.\footnote{In 
the case $(P,Q)=(3,1)$, there is also another solution $(x,y)=(1,1)$ with $y=1$, see Theorem~\ref{Tu31} for further details.}

Therefore, for some $n\in\mathbb{Z}$
$$\frac{v + u\sqrt{D}}2 = \pm \left(\frac{P + \sqrt{D}}2\right)^n$$
implying that
$$u = \pm U_n(P,Q)\qquad\text{and}\qquad v=\pm V_n(P,Q).$$
Since $u>0$ and $v>0$, we have
$$u = U_{|n|}(P,Q)\qquad\text{and}\qquad v=V_{|n|}(P,Q).$$
\end{proof}

\section{Quadratic Diophantine equations}\label{SecQuad}

While theorems in this section can be viewed as consequences of the fundamental results of Thue~\cite{Thue1909} and Siegel~\cite{Siegel1929},
we focus on computational aspects and give independent self-contained proofs of algorithmic nature.

\begin{Th}\label{Tqsys} 
A system of Diophantine equations
$$
\begin{cases} 
f_1(m,n) = k_1\\ 
f_2(m,n) = k_2 
\end{cases}
$$
where $f_1(m,n)$ and $f_2(m,n)$ are linearly independent homogeneous quadratic polynomials with integer coefficients and $(k_1, k_2)\ne (0,0)$ are integers, 
has a finite number of solutions.
\end{Th}

\begin{proof} Without loss of generality assume that $k_1\ne 0$.

If $m$ or $n$ is zero, then the other indeterminate may take only up to two values satisfying the quadratic equation 
$f_1(0,n)=k_1$ or $f_1(m,0)=k_1$. For the rest assume that both $m$ and $n$ are non-zero.

The given system implies an equation $g(m,n) = 0$ where $g(m,n)=k_2 f_1(m,n) - k_1 f_2(m,n)$ is a homogeneous quadratic polynomial.

If $g(m,n)=c mn$ for some coefficient $c\ne 0$, then there are no non-zero solutions.

If $g(m,n)\ne c mn$, without loss of generality we assume that the coefficient of $m^2$ in $g(m,n)$ is non-zero, implying that $g(m,n)=h(m/n) n^2$ 
where $h(x)$ is a quadratic polynomial. In this case, we have a quadratic equation $h(m/n)=0$ that determines up to two possible values of the ratio $m/n$.
For a fixed value of $m/n$, the equation $f_1(m,n) = k_1$ has only a finite number of integer solutions.

Therefore, the given system of Diophantine equations has only a finite number of solutions.
\end{proof}

\begin{Th}\label{Tinfsol}
A Diophantine equation $f(m,n)=k$, where $f(m,n)$ is a homogeneous polynomial of degree $4$ with integer coefficients and $k\ne 0$ is an integer, 
may have an infinite number of integer solutions only if $f(m,n)=c\cdot g(m,n)^2$, where $c$ is an integer and $g(m,n)$ is a polynomial of degree $2$
with integer coefficients.
\end{Th}

\begin{proof} Since the polynomial $f(m,n)$ is homogeneous, it can be represented as $f(m,n) = n^4 g(m/n)$ where $g(x)$ is a univariate polynomial with integer coefficients.
Consider a factorization of $g(x)$ over $\mathbb{Q}$ that 
implies a factorization $f(m,n)$ into irreducible homogeneous polynomials over $\mathbb{Z}$. 

If $f(m,n)$ has two linear factors $a_1 m + b_1 n$ and $a_2 m + b_2 n$ that are linearly independent, then the equation $f(m,n)=k$ implies 
a finite number of systems of the form:
$$
\begin{cases} 
a_1 m + b_1 n = d_1\\ 
a_2 m + b_2 n = d_2 
\end{cases}
$$
where $d_1, d_2$ run over the divisors\footnote{Unless specified otherwise, the divisors of an integer include both positive and negative divisors.} of $k$. 
Each of these systems has at most one solution, implying that the equation $f(m,n)=k$ has only a finite number of solutions. 
For the rest we assume that $f(m,n)$ has no two linearly independent linear factors. 
Then the polynomial $f(m,n)$ is not of the form $c\cdot g(m,n)^2$ only in the following four cases:
\begin{enumerate}
\item $f(m,n)$ is irreducible
\item $f(m,m) = g_2(m,n) h_2(m,n)$
\item $f(m,m) = g_1(m,n) h_3(m,n)$
\item $f(m,m) = g_1(m,n)^2 h_2(m,n)$
\end{enumerate}
where $g_i(m,n)$ and $h_i(m,n)$ are irreducible polynomials of degree $i$ with integer coefficients; and $g_2(m,n)$ and $h_2(m,n)$ are linearly independent.

In the case 1, the polynomial $f(m,n)$ is irreducible, implying that the Diophantine equation $f(m,n)=k$ has 
a finite number of solutions due to the Thue theorem~\cite{Thue1909}.

In the case 2, the equation $f(m,n)=k$ is equivalent to the union of the following systems of equations:
$$
\begin{cases} g_2(m,n) = d\\ 
h_2(m,n) = k/d \end{cases}
$$
where $d$ goes over the divisors of $k$. Each of these systems has only a finite number of solutions by Theorem~\ref{Tqsys}.
Therefore, the total number of solutions is also finite in this case.

In the case 3, the solutions to $f(m,n)=k$ must also satisfy a Thue equation $h_3(m,n) = d$ for some divisor $d$ of $k$.
Since for each such $d$ there is only a finite number of solutions, the total number of solutions to $f(m,n)=k$ is finite.

In the case 4, we let $g_2(m,n) = g_1(m,n)^2$ and proceed as in the case 2 above.
\end{proof}

\begin{Th}\label{Tcoeff} A Diophantine equation $f(m,n)=k$, where $k\ne 0$ and 
$$
f(m,n) = A (D^2 m^4 + E^2 n^4) + 2 B m n (D m^2 + E n^2) + C m^2 n^2
$$
with integer coefficients $A, B, C, D, E$,
may have an infinite number of integer solutions only if $2 A^2 D E + B^2 - A C = 0.$
\end{Th}

\begin{proof} According to Theorem~\ref{Tinfsol}, the equation $f(m,n)=k$ may have an infinite number of solutions
only if 
$$f(m,n)=c\cdot (u m^2 + v mn + w n^2)^2 = c (u^2 m^4 + w^2 n^4) + 2 c v m n (u m + w n) + c (2uw+v^2) m^2 n^2$$
for some integer coefficients $c,u,v,w$. 
In this case $A = cr^2$, $B=cvr$, $C=c(2uw+v^2)$, $D=u/r$, $E=w/r$ for some non-zero rational number $r$ and thus
$$
\begin{array}{lll}
2 A^2 D E + B^2 - A C & = & 2(cr^2)^2 (u/r) (w/r) + (cvr)^2 - (cr^2) c (v^2 + 2uw)\\
& = & c^2 r^2 (2uw + v^2 - (2uw + v^2)) = 0.
\end{array}
$$
\end{proof}

The theorem below represents a correct version of Corollary 6.3.8 in \cite{Cohen07} which appears seriously flawed.\footnote{
Applying Corollary 6.3.8~\cite{Cohen07} to the equation $-25x^2 + y^2 = -z^2$ with a particular solution $(x_0,y_0,z_0)=(1,0,5)$, all integer solutions 
to this equation are given by $(x,y,z)=\frac{u}{v}(-25s^2-t^2,-50st,-125s^2+5t^2)$, 
where the integers $s,t$ as well as $u,v$ are coprime and $v\mid 50$. 
However, it is easy seen that, under the last restriction, the solution $(x,y,z)=(1,3,4)$ cannot be obtained.}

\begin{Th}\label{ThABC}
Let $A, B, C$ be non-zero integers and let $(x_0, y_0, z_0)$ be a particular non-trivial integer solution of 
the Diophantine equation $Ax^2+By^2+Cz^2=0$ with $z_0 \ne 0$.
The general integer solution to the above equation is given by 
\begin{equation}\label{SolSys}
(x, y, z) = \frac{p}{q}\; ( P_x(m, n),\; P_y(m, n),\; P_z(m, n))
\end{equation}
where $m, n$ as well as $p, q$ are coprime integers with $q>0$ dividing $2 \lcm(A,B) Cz_0^2$, and
\begin{eqnarray*} 
P_x(m,n) & = & x_0 A m^2 + 2y_0 Bmn - x_0 B n^2,\\
P_y(m,n) & = & -y_0 A m^2 + 2x_0 A mn + y_0 B n^2,\\
P_z(m,n) & = & z_0 A m^2 + z_0 B n^2.
\end{eqnarray*}

\end{Th}

\begin{proof} The general solution~\eqref{SolSys} to the equation $Ax^2+By^2 = Cz^2$ is given in Corollary 6.3.6 in \cite{Cohen07}.
It remains to prove that $q$ divides $2 \lcm(A,B) Cz_0^2$.

The values of $x$, $y$, and $z$ in~\eqref{SolSys} are integer if and only if 
$q$ divides each of $P_x(m,n)$, $P_y(m,n)$, and $P_z(m,n)$. Then $q$ also divides
$$A(y_0 m - 2x_0 n)\cdot P_x(m,n) + (x_0 A m + 2y_0 B n)\cdot P_y(m,n) = 2BC z_0^2 n^3$$
as well as
$$(2 x_0 A m + y_0 B n)\cdot P_x(m,n) + (-2 y_0 B m + x_0 B n)\cdot P_y(m,n) = 2 A C z_0^2 m^3.$$
Therefore, $q$ divides
$$\gcd(2BC z_0^2 n^3, 2 A C z_0^2 m^3) = 2 \gcd(A,B) C z_0^2\gcd(\frac{B}{\gcd(A,B)} n^3, \frac{A}{\gcd(A,B)}m^3)$$
which in turn divides 
$$2 \gcd(A,B) C z_0^2 \frac{A}{\gcd(A,B)} \frac{B}{\gcd(A,B)} \gcd( n^3, m^3) = 2\lcm(A,B) C z_0^2.$$
\end{proof}

\begin{Th}\label{Txyz}
A system of Diophantine equations
\beq\label{Lxyzsys}
\begin{cases} 
a_1 x^2 + b_1 y^2 + c_1 z^2 = d_1\\
a_2 x^2 + b_2 y^2 + c_2 z^2 = d_2
\end{cases}
\eeq
where $a_i, b_i, c_i, d_i$ ($i=1,2$) are integers and in the matrix of coefficients
$$
M = \begin{bmatrix}
a_1 & b_1 & c_1 & d_1\\
a_2 & b_2 & c_2 & d_2
\end{bmatrix}
$$
every $2\times 2$ minor is non-zero, has a finite number of solutions.
\end{Th}

\begin{proof} It is clear that $(d_1,d_2)\ne(0,0)$. Without loss of generality we assume that $d_1\ne 0$.

If $x=0$ then the resulting system has a finite number of solutions due to Theorem~\ref{Tqsys}.
For the rest of the proof assume that $x\ne 0$.

Multiplying the first equation of \eqref{Lxyzsys} by $d_2$ and subtracting the second equation multiplied by $d_1$, we get a homogeneous Diophantine equation:
\beq\label{Lxyzhom}
A x^2 + B y^2 + C z^2 = 0
\eeq
where $A = a_1 d_2 - a_2 d_1$, $B = b_1 d_2 - b_2 d_1$, and $C = c_1 d_2 - c_2 d_1$ 
represent minors of the matrix $M$ and thus are non-zero.\footnote{From 
computational perspective, it is worth to divide each of the coefficients $A$, $B$, $C$ by their greatest common divisor $\gcd(A,B,C)$ before solving equation \eqref{Lxyzhom}.}

Let $(x_0,y_0,z_0)$ be a particular solution to equation \eqref{Lxyzhom} 
(for general methods of finding particular solutions to quadratic homogeneous equations in three variables, see \cite{Cremona03,Cohen07})\footnote{PARI/GP
provides a function \emph{bnfisnorm}, which can be used to compute a particular solution to a quadratic homogeneous equation 
in three indeterminants.}
with $z_0\ne 0$, i.e.,
\beq\label{Lxyz0}
A x_0^2 + B y_0^2 + C z_0^2 = 0,\qquad z_0\ne 0
\eeq
then by Theorem~\ref{ThABC} the general solution to equation \eqref{Lxyzhom} is given by \eqref{SolSys}, where 
$\gcd(m,n)=1$, $\gcd(p,q)=1$, and $q>0$ divides $2\lcm(A,B)Cz_0^2$.

Plugging the general solution \eqref{SolSys} into 
the first equation of \eqref{Lxyzsys}, we have
\beq\label{Lf12}
f(m,n) = d_1 \left(\frac{q}{p}\right)^2
\eeq
where
\begin{eqnarray*}
f(m,n) & \eqdef & (a_1 x_0^2 + b_1 y_0^2 + c_1 z_0^2)(A^2 m^4 + B^2 n^4)\\ 
& + & 4 x_0 y_0 (a_1 B - b_1 A) (A m^2 - B n^2) mn\\
& + & (4 a_1 y_0^2 B^2 - 2 A B (a_1 x_0^2 + b_1 y_0^2 - c_1 z_0^2) + 4 b_1 x_0^2 A^2) m^2 n^2.
\end{eqnarray*}
We remark that equation \eqref{Lf12} represents a collection of Diophantine equations indexed by coprime integers $p$ and $q$
ranging over the divisors of respectively $\sqrt{d_1 / d'_1}$ and $2\lcm(A,B)Cz_0^2$, where $d'_1$ is the square-free core of $d_1$.\footnote{Since 
the number of resulting equations is proportional to the number of divisors of $2\lcm(A,B)Cz_0^2$, 
from computational perspective it may be worth to choose a particular solution $(x_0,y_0,z_0)$ and re-label the variables (if needed) such that 
the number of divisors of $2\lcm(A,B)Cz_0^2$ is minimized.}

Theorem~\ref{Tcoeff} implies that equation \eqref{Lf12} may have infinitely many solutions only if 
\begin{align*}
&- 2 (a_1 x_0^2 + b_1 y_0^2 + c_1 z_0^2)^2 A B + (2 x_0 y_0 (a_1 B - b_1 A))^2 \\
&- (a_1 x_0^2 + b_1 y_0^2 + c_1 z_0^2) (4 a_1 y_0^2 B^2 - 2 A B (a_1 x_0^2 + b_1 y_0^2 - c_1 z_0^2) + 4 b_1 x_0^2 A^2) \\
= & -4 (A b_1 x_0^2 + B b_1 y_0^2 + B c_1 z_0^2) (A a_1 x_0^2 + B a_1 y_0^2 + A c_1 z_0^2) \\
= & - 4 z_0^4 (B c_1 - C b_1) (A c_1 - C a_1)\\
= & -4 z_0^4 d_1^2 (b_1 c_2 - b_2 c_1) (a_1 c_2 - a_2 c_1)
\end{align*}
is zero. However, each factor in this product is non-zero. 
Therefore, equation \eqref{Lf12} has only a finite number of integer solutions. 
\end{proof}

\section{Intersection of two Lucas sequences}\label{SecInter}

\subsection{Intersection of two $U$-sequences}\label{SecUU}

If an integer $x$ belongs to the intersection $\{ U_n(P_1,Q_1)\}$ and $\{ U_n(P_2,Q_2)\}$, where $|Q_1|=|Q_2|=1$ and $(P_1,Q_1)\ne (P_2,Q_2)$, 
then \eqref{MainEq} implies that 
$x$ satisfies a system of Diophantine equations:
\begin{equation}\label{MainSys}
\begin{cases} 
y^2 - D_1 x^2 = \pm 4\\
z^2 - D_2 x^2 = \pm 4
\end{cases}
\end{equation}
where $D_1 = P_1^2 - 4Q_1$ and $D_2 = P_2^2 - 4Q_2$. The matrix of its coefficients
$$
\begin{bmatrix}
-D_1 & 1 & 0 & \pm 4\\
-D_2 & 0 & 1 & \pm 4
\end{bmatrix}
$$
has all $2\times 2$ minors non-zero unless $D_1=D_2$ and the signs of 4's are the same.

If $D_1\ne D_2$, then system \eqref{MainSys} 
has only a finite number of solutions by Theorem~\ref{Txyz}. 
We note that it always has 
solutions $(x,y,z)=(0,2,2)$ and $(x,y,z)=(1,P_1,P_2)$ corresponding 
to common elements $0=U_0(P_1,Q_1)=U_0(P_2,Q_2)$ and 
$1=U_1(P_1,Q_1)=U_1(P_2,Q_2)$.

The case of $D_1 = D_2$ is addressed by the following theorem.

\begin{Th}\label{Tu31}
For distinct non-degenerate sequences $\{ U_n(P_1,Q_1) \}$ and $\{ U_n(P_2,Q_2) \}$ with $P_i>0$ and $|Q_i|=1$ ($i=1,2$), the following conditions are
equivalent:
\begin{itemize}
\item their intersection is infinite;
\item their discriminants are equal;
\item these sequences are Fibonacci numbers $\{ U_n(1,-1) \}$ and their bisection $\{ U_n(3,1) \} = \{ U_{2n}(1,-1) \}$.
\end{itemize}
\end{Th}

\begin{proof} As discussed above, the intersection of $\{ U_n(P_1,Q_1) \}$ and $\{ U_n(P_2,Q_2) \}$ can be infinite only if their discriminants are equal, i.e.,
$P_1^2 - 4Q_1 = P_2^2 - 4Q_2$.

If $Q_1=Q_2$ then $P_1 = P_2$, which is not the case for distinct sequences.

If $Q_1=-Q_2$ then $P_1^2 - P_2^2 = \pm 8$, implying that $(P_1,P_2)=(3,1)$ or $(P_1,P_2)=(1,3)$ and thus the sequences are $\{ U_n(1,-1) \}$ and $\{ U_n(3,1) \}$. 
Since $U_n(3,1)=U_{2n}(1,-1)$ for every integer $n$, $\{ U_n(3,1) \}$ represents a subsequence of $U_n(1,-1)$, implying that they have an infinite intersection:
$$
\{ U_n(1,-1) \} \cap \{ U_n(3,1) \} = \{ U_n(3,1) \}.
$$
\end{proof}

\subsection{Intersection of a $U$-sequence and a $V$-sequence}

If an integer $x$ belongs to the intersection $\{ U_n(P_1,Q_1)\}$ and $\{ V_n(P_2,Q_2)\}$, where $|Q_1|=|Q_2|=1$, then \eqref{MainEq} implies that 
$x$ satisfies a system of Diophantine equations:
\begin{equation}\label{UVSys}
\begin{cases} 
y^2 - D_1 x^2 = \pm 4\\
x^2 - D_2 z^2 = \pm 4
\end{cases}
\end{equation}
where $D_1 = P_1^2 - 4Q_1$ and $D_2 = P_2^2 - 4Q_2$. The matrix of its coefficients
$$
\begin{bmatrix}
-D_1 & 1 & 0 & \pm 4\\
1 & 0 & -D_2 & \pm 4
\end{bmatrix}
$$
has all $2\times 2$ minors non-zero unless $D_1=1$, which is not possible.
By Theorem~\ref{Txyz}, system \eqref{UVSys} has only a finite number of solutions.

\subsection{Intersection of two $V$-sequences}\label{SecPell}

If an integer $x$ belongs to the intersection $\{ V_n(P_1,Q_1)\}$ and $\{ V_n(P_2,Q_2)\}$, where $|Q_1|=|Q_2|=1$ and $(P_1,Q_1)\ne (P_2,Q_2)$, then \eqref{MainEq} implies that 
$x$ satisfies a system of Diophantine equations:
\begin{equation}\label{VSys}
\begin{cases} 
x^2 - D_1 y^2 = \pm 4\\
x^2 - D_2 z^2 = \pm 4
\end{cases}
\end{equation}
where $D_1 = P_1^2 - 4Q_1$ and $D_2 = P_2^2 - 4Q_2$. The matrix of its coefficients
$$
\begin{bmatrix}
1 & -D_1 & 0 & \pm 4\\
1 & 0 & -D_2 & \pm 4
\end{bmatrix}
$$
has all $2\times 2$ minors non-zero unless the signs of 4's in the r.h.s. are the same.
In the case of distinct signs in the r.h.s. of system \eqref{VSys}, it has a finite number of solutions by Theorem~\ref{Txyz}. 
Below we address the case of the same signs, 
which leads to the following system of equation:
\begin{equation}\label{VSys2}
\begin{cases} 
x^2 - D_1 y^2 = \pm 4\\
D_1 y^2 - D_2 z^2 = 0.
\end{cases}
\end{equation}

The second equation of system \eqref{VSys2} always has a solution $y=z=0$ corresponding to the common element $2$ of any two $V$-sequences.
On the other hand, the system has nonzero solutions only if the product $D_1 D_2$ is a square. 

\begin{Th} The intersection of distinct non-degenerate sequences 
$V_n(P_1,Q_1)$ and $V_n(P_2,Q_2)$ with $P_i>0$ and $|Q_i|=1$ ($i=1,2$) is infinite if and only if the product of 
their discriminants is a square. Moreover, in this case the intersection 
itself represents a Lucas sequence $V_n(P',Q')$ for some $P'>0$ and $|Q'|=1$.
\end{Th}

\begin{proof} As shown above, the intersection of $V_n(P_1,Q_1)$ and $V_n(P_2,Q_2)$ can be infinite only if the product of 
their discriminants $D_1=P_1^2 - 4Q_1$ and $D_2=P_2^2-4Q_2$ is a square. To prove the converse, assume that $D_1D_2$ is a square,
i.e., $D_1 = c d_1^2$ and $D_2 = c d_2^2$ for some integers $d_1$, $d_2$ with $\gcd(d_1,d_2)=1$, and non-square integer $c$. 

We first note that in this case, system \eqref{VSys} with opposite signs in the r.h.s. does not have a solution. 
Indeed, it implies that
$$x^4-16=(x^2 + 4)(x^2 - 4) = D_1 D_2 y^2 z^2 = (c d_1 d_2 y z)^2$$
and hence
$$(x^2 - c d_1 d_2 y z)(x^2 + c d_1 d_2 y z) = 16$$
whose the only solutions $x = \pm 2$ do not satisfy the original system \eqref{VSys}.
Therefore, the intersection of $\{ V_n(P_1,Q_1) \}$ and $\{ V_n(P_2,Q_2) \}$ necessarily represents the solutions to \eqref{VSys} with the same signs in the r.h.s.,
which is equivalent to \eqref{VSys2}.

The general solution to the second equation of \eqref{VSys2} is $(y,z)=\left(\pm d_2 t,\pm d_1 t\right)$ for $t\in\mathbb{Z}$, where the signs are independent of each other.
Plugging it into the first equation of \eqref{VSys2}, we get a Pellian equation:
\begin{equation}\label{EqPelld1d2}
x^2 - D t^2 = \pm 4
\end{equation}
where $D=\frac{D_1 D_2}{\gcd(D_1,D_2)} = c d_1^2 d_2^2.$

Notice that both $D_1$ and $D_2$ are congruent to $0$ or $1$ modulo $4$, and so is $D$. Therefore, it has infinitely many solutions $(x,t)=(x_k,t_k)$ of the form:
$$\frac{x_k + t_k\sqrt{D}}2 = \pm \left(\frac{x_1 + t_1\sqrt{D}}2\right)^k,\qquad k\in\mathbb{Z},$$
where $(x_1,t_1)$ is a solution with the smallest positive $t_1$ (see \cite{Cohen07}, Proposition 6.3.16). 
It follows that the sequence $x_k$ is simply $V_k(t_1,\mathrm{sign}(x_1^2 - D t_1^2))$.

Therefore, the intersection of Lucas sequences $V_n(P_1,Q_1)$ and $V_n(P_2,Q_2)$ forms 
Lucas sequence $V_k(t_1,\mathrm{sign}(x_1^2 - D t_1^2))$.

\end{proof}

\section{Computing intersections}\label{SecComp}

We note that the theorems in Section~\ref{SecBasic} have constructive proofs that allow straightforward implementation.
A bottleneck of our approach in practice is Thue equations that are somewhat hard to solve.
Nevertheless, there exist algorithms for solving Thue equations~\cite{Tzanakis1989,Bilu1996}, while
PARI/GP computer algebra system~\cite{PARI} provides an implementation of such an algorithm.
Below we present some numerical results obtained with the help of PARI/GP.

\subsection{Intersection of Fibonacci and Pell numbers}

\begin{Th} The only numbers that are both Fibonacci and Pell are $0$, $1$, $2$, and $5$.
\end{Th}

\begin{proof} To find the intersection of Fibonacci numbers $U_n(1,-1)$ and Pell numbers $U_n(2,-1)$,
we first compute their discriminants $D_1=5$ and $D_2=8$ and plug them into system \eqref{MainSys}:
\begin{equation}\label{FPSys}
\begin{cases} 
y^2 - 5 x^2 = \pm 4\\
z^2 - 8 x^2 = \pm 4
\end{cases}
\end{equation}
As it was shown in Section~\ref{SecUU}, this system has only a finite number of solutions.
We will find them all, following the outline of the proof of Theorem~\ref{Txyz}. First we obtain a homogeneous equation 
(depending on the choice of signs of 4's in the r.h.s.), resulting in two possible systems of equations:
$$\begin{cases}
y^2 - 5 x^2 = \pm 4\\
3x^2 + y^2 - z^2 = 0
\end{cases}\qquad
\begin{cases}
y^2 - 5 x^2 = \pm 4\\
- 13x^2 + y^2 + z^2 = 0
\end{cases}
$$
that we solve below.

Using a particular solution $(x_0,y_0,z_0)=(0,1,1)$ to the equation $3x^2 + y^2 - z^2 = 0$, from \eqref{SolSys} we obtain a general solution:
$$
\begin{cases}
x = 2mn\frac{p}{q}\\
y = (3m^2 - n^2)\frac{p}{q}\\
z = (3m^2+n^2) \frac{p}{q}
\end{cases}
$$
where the integers $m,n$ as well as $p,q$ are coprime, and $q>0$ divides 6.
Plugging these expressions into the first equation of \eqref{FPSys}, we obtain equations
$$9m^4 - 26m^2n^2 + n^4 = \pm \left(\frac{2q}{p}\right)^2$$
where $p$ divides $2$. Since the polynomial in the l.h.s. is irreducible, these equations are Thue equations.

According to PARI/GP, they have the following coprime solutions $(m,n)$ 
and corresponding triples $(x,y,z)$ (up to signs):
\begin{center}
\begin{tabular}{|l||l|l|}
\hline
$(p,q)$ & $(|m|,|n|)$ & $(|x|,|y|,|z|)$ \\
\hline\hline
$(2, 1)$ & $(0, 1)$ & $(0,2,2)$ \\
\hline
$(1, 2)$ & $(0, 2)$, $(1,1)$, $(1,5)$ & $(0,2,2)$, $(1,1,2)$, $(5,11,14)$ \\
\hline
$(2, 3)$ & $(1, 0)$ & $(0,2,2)$ \\
\hline
$(1, 6)$ & $(2, 0)$, $(1,3)$ & $(0,2,2)$, $(1,1,2)$, $(5,11,14)$ \\
\hline
\end{tabular}
\end{center}

Using a particular solution $(x_0,y_0,z_0)=(1,2,3)$ to the equation $-13x^2 + y^2 + z^2 = 0$, from \eqref{SolSys} we obtain a general solution:
$$\begin{cases}
x = (13m^2 - 4mn + n^2)\frac{p}{q}\\
y = (26m^2 - 26mn + 2 n^2)\frac{p}{q}\\
z = (39 m^2 - 3 n^2) \frac{p}{q}
\end{cases}
$$
where the integers $m,n$ as well as $p,q$ are coprime, and $q>0$ divides 234.
Plugging these expressions into the first equation of \eqref{FPSys}, we obtain equations
$$ - 169m^4 - 832m^3n + 570m^2n^2 - 64mn^3 - n^4 = \pm \left(\frac{2q}{p}\right)^2$$
where $p$ divides $2$. Since the polynomial in the l.h.s. is irreducible, these equations are Thue equations.
According to PARI/GP, they have the following coprime solutions $(m,n)$ (up to simultaneous change of signs) 
and corresponding triples $(|x|,|y|,|z|)$:
\begin{center}
\begin{tabular}{|l||l|l|}
\hline
$(p,q)$ & $(m,n)$ & $(|x|,|y|,|z|)$ \\
\hline\hline
$(2, 1)$ & $(0, 1)$ & $(2,4,6)$ \\
\hline
$(1, 2)$ & $(0, 2)$ & $(2,4,6)$ \\
\hline
$(2, 9)$ & $(0, 3)$, $(1, 2)$ & $(2,4,6)$ \\
\hline
$(2, 13)$ & $(1, 0)$ & $(2,4,6)$ \\
\hline
$(1, 18)$ & $(0, 6)$, $(1, -1)$, $(1, 5)$, $(2, 4)$ & $(1,3,2)$, $(2,4,6)$ \\
\hline
$(1, 26)$ & $(2, 0)$ & $(2,4,6)$ \\
\hline
$(2, 117)$ & $(2, 13)$, $(3, 0)$ & $(2,4,6)$ \\
\hline
$(1, 234)$ & $(1, -13)$, $(5, 13)$, $(4, 26)$, $(6, 0)$ & $(1,3,2)$, $(2,4,6)$ \\
\hline
\end{tabular}
\end{center}

\end{proof}

\subsection{Intersections of some Lucas sequences}

In this section we lists some interesting numerical results and pose related open questions.
As before, we consider only on non-degenerate sequence and in addition exclude the sequences with $(P,Q)=(3,1)$ corresponding to bisections of those with $(P,Q)=(1,-1)$.

We call the intersection of sequences $U_n(P_1,Q_1)$ and $U_n(P_2,Q_2)$ \emph{non-trivial} if it contains an integer $m>1$ different from $P_1$ and $P_2$.
Examples of non-trivial finite intersections are 
$$
\begin{array}{lllll}
\{ U_n(3,-1) \} & \cap & \{ U_n(19,1) \} & = & \{ 0, 1, 360 \}, \\
\{ U_n(8,-1) \} & \cap & \{ U_n(23,1) \} & = & \{ 0, 1, 528 \}, \\
\{ U_n(3,-1) \} & \cap & \{ U_n(6,1) \} & = & \{ 0, 1, 1189 \}.
\end{array}
$$
For an integer $m>1$, we can find all values of $P>0$ such that $m\in\{ U_n(P,Q)\}$ for some $Q=\pm 1$ by solving the equation (followed from Theorem~\ref{Trecog})
$$y^2 - m^2 P^2 = \pm 4(m^2\pm 1).$$
In particular, the only integers $m\leq 10^9$ that belong to two distinct sequences 
$U_n(P_1,Q_1)$ and $U_n(P_2,Q_2)$ for some $P_1\ne m$, $P_2 \ne m$, and $|Q_1|=|Q_2|=1$, are $5$, $360$, $528$, and $1189$. 
Are there any integers $m>10^9$ with the analogous property?

Another notable intersection is
$$\{ U_n(1,-1) \} \cap \{ U_n(2,-1) \} = \{ 0, 1, 2, 5 \}$$
which is the only intersection of size four that we are aware about. Does there exist any other?

Similarly, we call the intersection of $U_n(P_1,Q_1)$ and $V_n(P_2,Q_2)$ \emph{non-trivial} if it contains 
an integer $m>2$ different from $P_1$, and $P_2$. 
Examples of non-trivial intersections are
$$\{ U_n(1,-1) \} \cap \{ V_n(2,-1) \} = \{ U_n(1,-1) \} \cap \{ V_n(6,1) \} = \{ 2, 34 \},$$
$$\{ U_n(9,-1) \} \cap \{ V_n(2,-1) \} = \{ 82 \}.$$
In particular, the only integers $m$, $2 < m\leq 10^9$, 
that belong the intersection $U_n(P_1,Q_1)$ and $V_n(P_2,Q_2)$ for some $P_1 \ne m$, $P_2\ne m$, and $|Q_1|=|Q_2|=1$,
are $3$, $29$, $34$, and $82$. 
Are there any integers $m>10^9$ with the analogous property?

The intersections
$$\{ U_n(1,-1) \} \cap \{ V_n(1,-1) \} = \{ 1, 2, 3 \},$$
$$\{ U_n(2,-1) \} \cap \{ V_n(1,-1) \} = \{ 1, 2, 29 \}$$
have size $3$. Are there any other?

Again, we call the intersection of $V_n(P_1,Q_1)$ and $V_n(P_2,Q_2)$ \emph{non-trivial} if it contains an 
integer different from $m>2$, $P_1$, and $P_2$. Examples of non-trivial finite intersections are
$$\{ V_n(1,-1) \} \cap \{ V_n(3,-1) \} = \{ 2, 3, 11 \},$$
$$\{ V_n(2,-1) \} \cap \{ V_n(4,1) \} = \{ 2, 14 \}.$$
Are there any other non-trivial finite intersections?

Table~\ref{TinfVV1} 
gives all parameters $(P_1,Q_1)$ and $(P_2,Q_2)$ such that $P_1<P_2\leq 1000$, $|Q_1|=|Q_2|=1$ and 
$\{ V_n(P_1,Q_1)\} \cap \{ V_n(P_2,Q_2) \}$ is infinite.

\begin{scriptsize}
\begin{table}[t]
\begin{center}
\begin{scriptsize}
\begin{tabular}{|l|l|l|l|l|l|l|}
\hline
$(P_1,Q_1)$ & $(P_2,Q_2)$ & $(P,Q)$  & &  $(P_1,Q_1)$ & $(P_2,Q_2)$ & $(P,Q)$ \\
\hline
(1,-1) & (3,1) & (3,1) & &                       (5,1) & (23,1) & (23,1) \\
(1,-1) & (4,-1) & (4,-1) & &                     (5,1) & (110,1) & (110,1) \\
(1,-1) & (7,1) & (7,1) & &                       (5,1) & (527,1) & (527,1) \\
(1,-1) & (11,-1) & (11,-1) & &                   (5,-1) & (27,1) & (27,1) \\
(1,-1) & (18,1) & (18,1) & &                     (5,-1) & (140,-1) & (140,-1) \\
(1,-1) & (29,-1) & (29,-1) & &                   (5,-1) & (727,1) & (727,1) \\
(1,-1) & (47,1) & (47,1) & &                     (6,1) & (14,-1) & (198,1) \\
(1,-1) & (76,-1) & (76,-1) & &                   (6,1) & (34,1) & (34,1) \\
(1,-1) & (123,1) & (123,1) & &                   (6,1) & (82,-1) & (6726,1) \\
(1,-1) & (199,-1) & (199,-1) & &                 (6,1) & (198,1) & (198,1) \\
(1,-1) & (322,1) & (322,1) & &                   (6,1) & (478,-1) & (228486,1) \\
(1,-1) & (521,-1) & (521,-1) & &                 (6,-1) & (38,1) & (38,1) \\
(1,-1) & (843,1) & (843,1) & &                   (6,-1) & (234,-1) & (234,-1) \\
(2,-1) & (6,1) & (6,1) & &                       (7,1) & (11,-1) & (15127,1) \\
(2,-1) & (14,-1) & (14,-1) & &                   (7,1) & (18,1) & (322,1) \\
(2,-1) & (34,1) & (34,1) & &                     (7,1) & (29,-1) & (710647,1) \\
(2,-1) & (82,-1) & (82,-1) & &                   (7,1) & (47,1) & (47,1) \\
(2,-1) & (198,1) & (198,1) & &                   (7,1) & (76,-1) & (33385282,1) \\
(2,-1) & (478,-1) & (478,-1) & &                 (7,1) & (123,1) & (15127,1) \\
(3,1) & (4,-1) & (18,1) & &                      (7,1) & (199,-1) & (1568397607,1) \\
(3,1) & (7,1) & (7,1) & &                        (7,1) & (322,1) & (322,1) \\
(3,1) & (11,-1) & (123,1) & &                    (7,1) & (521,-1) & (73681302247,1) \\
(3,1) & (18,1) & (18,1) & &                      (7,1) & (843,1) & (710647,1) \\
(3,1) & (29,-1) & (843,1) & &                    (7,-1) & (51,1) & (51,1) \\
(3,1) & (47,1) & (47,1) & &                      (7,-1) & (364,-1) & (364,-1) \\
(3,1) & (76,-1) & (5778,1) & &                   (8,1) & (62,1) & (62,1) \\
(3,1) & (123,1) & (123,1) & &                    (8,1) & (488,1) & (488,1) \\
(3,1) & (199,-1) & (39603,1) & &                 (8,-1) & (66,1) & (66,1) \\
(3,1) & (322,1) & (322,1) & &                    (8,-1) & (536,-1) & (536,-1) \\
(3,1) & (521,-1) & (271443,1) & &                (9,1) & (79,1) & (79,1) \\
(3,1) & (843,1) & (843,1) & &                    (9,1) & (702,1) & (702,1) \\
(3,-1) & (11,1) & (11,1) & &                     (9,-1) & (83,1) & (83,1) \\
(3,-1) & (36,-1) & (36,-1) & &                   (9,-1) & (756,-1) & (756,-1) \\
(3,-1) & (119,1) & (119,1) & &                   (10,1) & (98,1) & (98,1) \\
(3,-1) & (393,-1) & (393,-1) & &                 (10,1) & (970,1) & (970,1) \\
(4,1) & (14,1) & (14,1) & &                      (10,-1) & (102,1) & (102,1) \\
(4,1) & (52,1) & (52,1) & &                      (11,1) & (36,-1) & (1298,1) \\
(4,1) & (194,1) & (194,1) & &                    (11,1) & (119,1) & (119,1) \\
(4,1) & (724,1) & (724,1) & &                    (11,1) & (393,-1) & (154451,1) \\
(4,-1) & (7,1) & (322,1) & &                     (11,-1) & (18,1) & (1860498,1) \\
(4,-1) & (11,-1) & (1364,-1) & &                 (11,-1) & (29,-1) & (20633239,-1) \\
(4,-1) & (18,1) & (18,1) & &                     (11,-1) & (47,1) & (228826127,1) \\
(4,-1) & (29,-1) & (24476,-1) & &                (11,-1) & (76,-1) & (2537720636,-1) \\
(4,-1) & (47,1) & (103682,1) & &                 (11,-1) & (123,1) & (123,1) \\
(4,-1) & (76,-1) & (76,-1) & &                   (11,-1) & (199,-1) & (312119004989,-1) \\
(4,-1) & (123,1) & (1860498,1) & &               (11,-1) & (322,1) & (3461452808002,1) \\
(4,-1) & (199,-1) & (7881196,-1) & &             (11,-1) & (521,-1) & (38388099893011,-1) \\
(4,-1) & (322,1) & (322,1) & &                   (11,-1) & (843,1) & (425730551631123,1) \\
(4,-1) & (521,-1) & (141422324,-1) & &           (12,1) & (142,1) & (142,1) \\
(4,-1) & (843,1) & (599074578,1) & &             (12,-1) & (146,1) & (146,1) \\
\hline
\end{tabular}
\caption{Parameters such that $\{ V_n(P_1,Q_1) \} \cap \{ V_n(P_2,Q_2) \} = \{ V_n(P,Q) \}$.}
\label{TinfVV1}
\end{scriptsize}
\end{center}
\end{table}

\setcounter{table}{0}

\begin{table}[t]
\begin{center}
\begin{scriptsize}
\begin{tabular}{|l|l|l|l|l|l|l|}
\hline
$(P_1,Q_1)$ & $(P_2,Q_2)$ & $(P,Q)$  & &  $(P_1,Q_1)$ & $(P_2,Q_2)$ & $(P,Q)$ \\
\hline
(13,1) & (167,1) & (167,1) & &                   (29,-1) & (199,-1) & (12360848946698171,-1) \\
(13,-1) & (171,1) & (171,1) & &                  (29,-1) & (322,1) & (358890350005878082,1) \\
(14,1) & (52,1) & (2702,1) & &                   (29,-1) & (521,-1) & (10420180999117162549,-1) \\
(14,1) & (194,1) & (194,1) & &                   (29,-1) & (843,1) & (843,1) \\
(14,1) & (724,1) & (524174,1) & &                (30,1) & (898,1) & (898,1) \\
(14,-1) & (34,1) & (39202,1) & &                 (30,-1) & (902,1) & (902,1) \\
(14,-1) & (82,-1) & (551614,-1) & &              (31,1) & (959,1) & (959,1) \\
(14,-1) & (198,1) & (198,1) & &                  (31,-1) & (963,1) & (963,1) \\
(14,-1) & (478,-1) & (109216786,-1) & &          (34,1) & (82,-1) & (45239074,1) \\
(15,1) & (223,1) & (223,1) & &                   (34,1) & (198,1) & (39202,1) \\
(15,-1) & (227,1) & (227,1) & &                  (34,1) & (478,-1) & (52205852194,1) \\
(16,1) & (254,1) & (254,1) & &                   (36,-1) & (119,1) & (1684802,1) \\
(16,-1) & (258,1) & (258,1) & &                  (36,-1) & (393,-1) & (60699636,-1) \\
(17,1) & (287,1) & (287,1) & &                   (38,1) & (234,-1) & (54758,1) \\
(17,-1) & (291,1) & (291,1) & &                  (47,1) & (76,-1) & (1114577054219522,1) \\
(18,1) & (29,-1) & (599074578,1) & &             (47,1) & (123,1) & (228826127,1) \\
(18,1) & (47,1) & (103682,1) & &                 (47,1) & (199,-1) & (2459871053643326447,1) \\
(18,1) & (76,-1) & (5778,1) & &                  (47,1) & (322,1) & (103682,1) \\
(18,1) & (123,1) & (1860498,1) & &               (47,1) & (521,-1) & (5428934300813767249007,1) \\
(18,1) & (199,-1) & (62113250390418,1) & &       (47,1) & (843,1) & (505019158607,1) \\
(18,1) & (322,1) & (322,1) & &                   (51,1) & (364,-1) & (132498,1) \\
(18,1) & (521,-1) & (20000273725560978,1) & &    (52,1) & (194,1) & (7300802,1) \\
(18,1) & (843,1) & (599074578,1) & &             (52,1) & (724,1) & (379501252,1) \\
(18,-1) & (326,1) & (326,1) & &                  (62,1) & (488,1) & (238142,1) \\
(19,1) & (359,1) & (359,1) & &                   (66,1) & (536,-1) & (287298,1) \\
(19,-1) & (363,1) & (363,1) & &                  (76,-1) & (123,1) & (6440026026380244498,1) \\
(20,1) & (398,1) & (398,1) & &                   (76,-1) & (199,-1) & (489526700523968661124,-1) \\
(20,-1) & (402,1) & (402,1) & &                  (76,-1) & (322,1) & (33385282,1) \\
(21,1) & (439,1) & (439,1) & &                   (76,-1) & (521,-1) & (2828485190904971853895196,-1) \\
(21,-1) & (443,1) & (443,1) & &                  (76,-1) & (843,1) & (215002084978043708894524818,1) \\
(22,1) & (482,1) & (482,1) & &                   (79,1) & (702,1) & (492802,1) \\
(22,-1) & (486,1) & (486,1) & &                  (82,-1) & (198,1) & (304278004998,1) \\
(23,1) & (110,1) & (12098,1) & &                 (82,-1) & (478,-1) & (24954506565518,-1) \\
(23,1) & (527,1) & (527,1) & &                   (83,1) & (756,-1) & (571538,1) \\
(23,-1) & (531,1) & (531,1) & &                  (98,1) & (970,1) & (940898,1) \\
(24,1) & (574,1) & (574,1) &&                    (110,1) & (527,1) & (146361602,1) \\
(24,-1) & (578,1) & (578,1) &&                   (119,1) & (393,-1) & (23855111399,1) \\
(25,1) & (623,1) & (623,1) & &                   (123,1) & (199,-1) & (97418273275323406890123,1) \\
(25,-1) & (627,1) & (627,1) &&                   (123,1) & (322,1) & (3461452808002,1) \\
(26,1) & (674,1) & (674,1) &&                    (123,1) & (521,-1) & (1473646213395791149646646123,1) \\
(26,-1) & (678,1) & (678,1) &&                   (123,1) & (843,1) & (425730551631123,1) \\
(27,1) & (140,-1) & (19602,1) &&                 (140,-1) & (727,1) & (384238402,1) \\
(27,1) & (727,1) & (727,1) &&                    (194,1) & (724,1) & (274758382274,1) \\
(27,-1) & (731,1) & (731,1) &&                   (198,1) & (478,-1) & (11928306344169798,1) \\
(28,1) & (782,1) & (782,1) &&                    (199,-1) & (322,1) & (3858055874062761829426214722,1) \\
(28,-1) & (786,1) & (786,1) &&                   (199,-1) & (521,-1) & (767772505664398093937756525279,-1) \\
(29,1) & (839,1) & (839,1) &&                    (199,-1) & (843,1) & (152790586683089283455442974745243,1) \\
(29,-1) & (47,1) & (505019158607,1) &&           (322,1) & (521,-1) & (400010949097364802732720796316482,1) \\
(29,-1) & (76,-1) & (14662949395604,-1) &&       (322,1) & (843,1) & (358890350005878082,1) \\
(29,-1) & (123,1) & (425730551631123,1) &&       (521,-1) & (843,1) & (108580172054362347934782139964888177403,1) \\
\hline
\end{tabular}
\caption{(con'd) Parameters such that $\{ V_n(P_1,Q_1) \} \cap \{ V_n(P_2,Q_2) \} = \{ V_n(P,Q) \}$.}
\label{TinfVV2}
\end{scriptsize}
\end{center}
\end{table}

\end{scriptsize}

\section{Lucas sequences with arbitrary initial terms}\label{SecArbit}

For fixed $P, Q\in\mathbb{Z}$, $|Q|=1$, consider a sequence that satisfy the same recurrent relation as $U_n(P,Q)$ but with initial terms $a,b\in\mathbb{Z}$:
$$U^{a,b}_0(P,Q)=a,\quad U^{a,b}_1(P,Q)=b,\quad U^{a,b}_{n+1}(P,Q) = P\cdot U^{a,b}_n(P,Q) - Q\cdot U^{a,b}_{n-1}(P,Q),\, n\geq 1.$$
It is easy to see that it can be expressed in terms of conventional Lucas sequences as follows:
$$U^{a,b}_n(P,Q) = (b-aP)\cdot U_n(P,Q) + a\cdot U_{n+1}(P,Q).$$
This representation allows us to define a counterpart sequence:
$$V^{a,b}_n(P,Q) \eqdef (b-aP)\cdot V_n(P,Q) + a\cdot V_{n+1}(P,Q) = U^{2b-aP,bP-2aQ}_n(P,Q).$$

It is easy to see that the generalized Lucas sequences satisfy an identity similar to \eqref{MainEq}:
$$V^{a,b}_n(P,Q)^2 - D\cdot U^{a,b}_n(P,Q)^2 = ((2b-aP)^2 - D a^2)Q^n$$
where $D=P^2-4Q$. However, there is no general analog of Theorem~\ref{Trecog} and hence the identity
\begin{equation}\label{EqPellAB}
y^2 - D\cdot x^2 = \pm ((2b-aP)^2 - D a^2)
\end{equation}
can be used only as a necessary condition for testing membership of $x$ in $U^{a,b}_n(P,Q)$. 
In particular, while $(x,y)=(8,34)$ is a solution to equation \eqref{EqPellAB} for $(P,Q,a,b)=(1,-1,1,15)$, 8 is not an element of $U^{1,15}_n(1,-1)$.

The quadratic equation \eqref{EqPellAB} allows one to compute the intersection of Lucas sequences with arbitrary initial terms and prove its finiteness similarly to the case of classic Lucas sequences.

\section*{Acknowledgements}

The author thanks Rustem Aydagulov for a number of insightful discussions.

\begin{small}
\bibliographystyle{plain} 
\bibliography{fib-pell.bib} 
\end{small}

\end{document}